\documentclass[11pt]{article}
\usepackage{amsfonts,amsmath,latexsym,color,epsfig}
\usepackage{epsfig}

\date{}

\title{A crossing lemma for multigraphs}

\author{
{\sl J\'anos Pach}\thanks{
Ecole Polytechnique F\'ed\'erale de Lausanne and
R\'enyi Institute, Hungarian Academy of Sciences,
P.O.Box 127 Budapest, 1364, Hungary; \texttt{pach@renyi.hu}. ;
\texttt{pach@cims.nyu.edu}.
Supported by Swiss National Science Foundation Grants
200021-165977 and 200020-162884 and Schloss Dagstuhl -- Leibniz Center for Informatics.}
\and
{\sl G\'eza T\'oth}\thanks{
R\'enyi Institute, Hungarian Academy of Sciences,
P.O.Box 127 Budapest, 1364, Hungary; \texttt{geza@renyi.hu}.
Supported by National Research, Development and Innovation Office, NKFIH,
K-111827 and Schloss Dagstuhl -- Leibniz Center for Informatics.}}

\begin{document}

\maketitle

\begin{abstract}
Let $G$ be a drawing of a graph with $n$ vertices and $e>4n$ edges, in which no two adjacent edges
cross and any pair of independent edges cross at most once. According to the celebrated Crossing Lemma
of Ajtai, Chv\'atal, Newborn, Szemer\'edi and Leighton, the number of crossings in $G$ is at least
$c{e^3\over n^2}$, for a suitable constant $c>0$. In a seminal paper, Sz\'ekely generalized this
result to multigraphs, establishing the lower bound $c{e^3\over mn^2}$, where $m$ denotes the maximum
multiplicity of an edge in $G$. We get rid of the dependence on $m$ by showing that, as in the
original Crossing Lemma, the number of crossings is at least $c'{e^3\over
  n^2}$ for some $c'>0$, provided that the ``lens'' enclosed by every pair of
parallel edges in $G$ contains at least one vertex. This settles a conjecture
of
Kaufmann.
\end{abstract}

\section{Introduction}

A {\em drawing} of a graph $G$ is a representation of $G$ in the plane such
that the vertices are represented by points, the edges are represented by
simple continuous arcs connecting the corresponding pair of points without
passing through any other point representing a vertex. In notation and
terminology we do not make any distinction between a vertex (edge) and the
point (resp., arc) representing it. Throughout this note we assume that any
pair of edges intersect in finitely many points and no three edges pass
through the same point. A common interior point of two edges at which the
first edge passes from one side of the second edge to the other, is called
a {\em crossing}.

A very ``successful concept for measuring non-planarity" of graphs is the
{\em crossing number} of $G$~\cite{Sz04}, which is defined as the minimum
number ${\rm cr}(G)$ of crossing points in any drawing of $G$ in the plane.
For many interesting variants of the crossing number, see  ~\cite{Sch13},
\cite{PT00b}. Computing ${\rm cr}(G)$ is an NP-hard problem~\cite{GJ83},
which is equivalent to the existential theory of reals~\cite{Sch10}.

The following statement, proved independently by Ajtai, Chv\'atal, Newborn,
Szemer\'edi~\cite{ACNS82} and Leighton~\cite{L83}, gives a lower bound on the
crossing number of a graph in terms of its number of vertices and number of
edges.

\medskip

\noindent {\bf Crossing Lemma.} \cite{ACNS82}, \cite{L83}
{\em   For any graph $G$ with $n$ vertices and $e>4n$ edges, we have
 $${\rm cr}(G) \geq \frac{1}{64}\frac{e^3}{n^2}.$$}

\medskip
Apart from the exact value of the constant, the order of magnitude of this
bound cannot be improved. This lemma has many important applications,
including simple proofs of the Szemer\'edi-Trotter theorem~\cite{SzT83} on the
maximum number of incidences between $n$ points and $n$ lines in the plane and
of the best known upper bound on the number of halving lines induced by $n$
points, due to Dey \cite{D98}.

The same problem was also considered for {\em multigraphs} $G$, in which two
vertices can be connected by several edges. As Sz\'ekely~\cite{Sz97} pointed
out, if the {\em multiplicity} of an edge is at most $m$, that is, any pair of
vertices of $G$ is connected by at most $m$ ``parallel'' edges, then the
minimum number of crossings between the edges satisfies
\begin{equation}\label{eq1}
{\rm cr}(G) \geq \frac{1}{64}\frac{e^3}{mn^2}
\end{equation}
when $e\ge 4mn$.
For $m=1$, this gives the Crossing Lemma, but as $m$ increases, the bound is
getting weaker. It is not hard to see that this inequality is also tight up to
a constant factor. Indeed, consider any (simple) graph with $n$ vertices and
roughly $e/m>4n$ edges such that it can be drawn with at
most $\frac{(e/m)^3}{n^2}$ crossings, and replace each edge by $m$ parallel
edges no pair of which share an interior point. The crossing number of the
resulting multigraph cannot exceed $\frac{(e/m)^3}{n^2}m^2=\frac{e^3}{mn^2}$.

\smallskip

It was suggested by Michael Kaufmann~\cite{Ka16} that the dependence on the
multiplicity might be eliminated if we restrict our attention to a special
class of drawings.

\medskip

\noindent {\bf Definition.}
{A {\em drawing} of a multigraph $G$ in the plane is called
{\em branching}, or a {\em branching topological multigraph}, if the following
conditions are satisfied.

(i) If two edges are parallel (have the same endpoints), then there is at
least one vertex in the interior and in the exterior of the simple closed
curve formed by their union.

(ii) If two edges share at least one endpoint, they cannot cross.

(iii) If two edges do not share an endpoint, they can have at most one crossing.
\smallskip

Given a multigraph $G$, its {\em branching crossing number} is the smallest
number ${\rm cr}_{\rm br}(G)$ of crossing points in any branching drawing of
$G$. If $G$ has no such drawing, set ${\rm cr}_{\rm br}(G)=\infty$.}
\medskip

According to this definition,  ${\rm cr}_{\rm br}(G)\ge{\rm cr}(G)$
for every graph or multigraph $G$, and if $G$ has no parallel edges, equality
holds.

The main aim of this note is to settle Kaufmann's conjecture.

\medskip
\noindent {\bf Theorem 1.}
{\em The branching crossing number of any multigraph $G$ with
$n$ vertices and $e>4n$ edges satisfies
${\rm cr}_{\rm br}(G) \geq c\frac{e^3}{n^2}$, for an absolute constant $c>10^{-7}$.}

\medskip

Unfortunately, the standard proofs of the Crossing Lemma by inductional or
probabilistic arguments break down in this case, because the property that a
drawing of $G$ is branching is not hereditary: it can be destroyed by deleting
vertices from $G$.

\smallskip

The bisection width of an {\em abstract} graph is usually defined as the
minimum number of edges whose deletion separates the graph into two parts
containing ``roughly the same'' number of vertices. In analogy to this, we
introduce the following new parameter of {\em branching topological}
multigraphs.

\medskip

\noindent {\bf Definition.}
{The {\em branching bisection width} ${\rm b}_{\rm br}(G)$ of a {\em
    branching topological multigraph} $G$ with $n$ vertices is the minimum
  number of edges whose removal
splits
$G$ into two {\em branching topological multigraphs}, $G_1$ and $G_2$, with no
  edge connecting them such that $|V(G_1)|, |V(G_2)|\ge n/5$.}

\medskip

A key element of the proof of Theorem 1 is the following statement
establishing a relationship between the branching bisection width and the
number of crossings of a branching topological multigraph.

\medskip
\noindent {\bf Theorem 2.} {\em Let $G$ be a branching topological multigraph
with $n$ vertices of degrees $d_1, d_2,\ldots, d_n$, and with $c(G)$
crossings. Then the branching bisection width of $G$ satisfies
$${\rm b}_{\rm br}(G)\le 22\sqrt{c(G)+\sum_{i=1}^nd_i^2+n}.$$}
\medskip

By definition, the number of crossings $c(G)$ between the edges of $G$ has to
be at least as large as the branching crossing number of the abstract
underlying multigraph of $G$.

\smallskip

To prove Theorem 1, we will use Theorem 2 recursively. Therefore, it is
crucially important that in the definition of ${\rm b}_{\rm br}(G)$, both
parts that $G$ is cut into should be branching topological multigraphs
themselves. If we are not careful, all vertices of $V(G)$ that lie in the
interior (or in the exterior) of a closed curve formed by two parallel edges
between $u,v\in G_1$, say, may end up in $G_2$. This would violate for $G_1$
condition (i) in the above definition of branching topological
multigraphs. That is why the proof of Theorem 2 is far more delicate than the
proof of the analogous statement for abstract graphs without multiple edges,
obtained in ~\cite{PSS96}.

\smallskip

For the proof of Theorem 1, we also need the following result.

\medskip
\noindent {\bf Theorem 3.} {\em Let $G$ be a branching topological multigraph
  with $n\ge 3$ vertices. Then the number of edges of $G$ satisfies $e(G)\le
  n(n-2)$, and this bound is tight.}

\vskip 0.3cm

\begin{figure}[h]
 \begin{center}
  \includegraphics[scale=0.5]{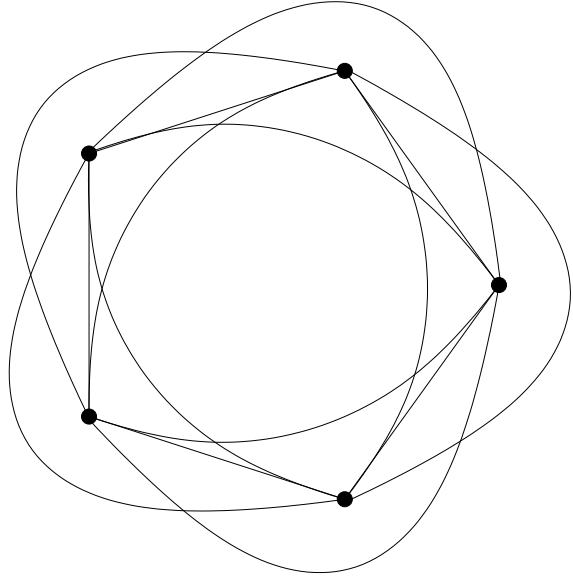}
 \end{center}
\caption{Theorem 3 is tight for every $n\ge 3$. Construction for $n=5$.}\label{otszog}
\end{figure}

\medskip

Our strategy for proving Theorem 1 is the following. Suppose, for a
contradiction, that a multigraph $G$ has a branching drawing in which the
number of crossings is smaller than what is required by the theorem. According
to Theorem 2, this implies that the branching bisection width of this drawing
is small. Thus, we can cut the drawing into two smaller branching topological
multigraphs, $G_1$ and $G_2$, by deleting relatively few edges. We repeat the
same procedure for $G_1$ and $G_2$, and continue recursively until the size of
every piece falls under a carefully chosen threshold. The total number of
edges removed during this procedure is small, so that the small components
altogether still contain a lot of edges. However, the number of edges in the
small components is bounded from above by Theorem 3, which leads to the
desired contradiction.

\medskip
\noindent{\bf Remarks.}
{\bf 1.}  Theorem 1 does not hold if we drop conditions (ii) and (iii) in the above definition,
that is, if we allow two edges to cross more than once. To see this, suppose
that $n$ is a multiple of 3 and consider a tripartite topological multigraph
$G$ with $V(G)=V_1\cup V_2\cup V_3$, where all points of $V_i$ belong to the
line $x=i$ and we have $|V_i|=n/3$ for $i=1,2,3$. Connect each point of $V_1$
to every point of $V_3$ by $n/3$ parallel edges: by one curve passing
between any two (cyclically) consecutive vertices of $V_2$. We can draw these curves in
such a way that any two edges cross at most twice, so that the
number of edges is $e=e(G)=(n/3)^3$ and the total number of crossings is at
most $2{e\choose 2}<(n/3)^6$. On the other hand, the
lower bound in Theorem~1 is $ce^3/n^2>(c/3^{9})n^7$, which is a
contradiction if $n$ is sufficiently large.

{\bf 2.}   In the definition of {\em branching topological multigraphs}, for symmetry
we assumed that the closed curve obtained by the concatenation of any pair of parallel edges in $G$
has at least one vertex in its interior and at least one vertex in its exterior; see condition (i).
It would have been sufficient to require that any such curve has at least one vertex in its
{\em interior}, that is, any lens enclosed by two parallel edges contains a vertex. Indeed, by
placing an isolated vertex $v$ far away from the rest of the drawing, we can achieve that there is
at least one vertex (namely, $v$) in the exterior of every lens, and apply Theorem 1 to the resulting
graph with $n+1$ vertices.

{\bf 3.}   Throughout this paper, we assume for simplicity that a multigraph does not have {\em loops},
that is, there are no edges whose endpoints are the same. It is easy to see that Theorem 1, with a slightly worse constant $c$,
also holds for topological multigraphs $G$ having loops, provided that condition (ii) in the definition of branching
topological multigraphs remains valid. In this case, one can argue that the total number of loops
cannot exceed $n$. Subdividing every loop by an additional vertex, we get rid of all loops, and then
we can apply Theorem 1 to the resulting multigraph of at most $2n$ vertices.

\smallskip
The rest of this note is organized as follows. In Section 2, we establish
Theorem 3. In Section 3, we apply Theorems 2 and 3 to deduce Theorem 1. The
proof of Theorem 2 is given in Section 4.

\section{The number of edges in branching topological multigraphs\\
---Proof of Theorem 3}

\noindent {\bf Lemma 2.1.}
{\em Let $G$ be a branching topological multigraph with $n\ge 3$ vertices and
  $e$ edges, in which no two edges cross each other.
Then $e\le 3n-6$.}

\smallskip

\noindent {\bf Proof.}
We can suppose without loss of generality that $G$ is connected. Otherwise, we
can achieve this by adding some edges of multiplicity $1$, without violating
conditions (i)-(iii) required for a drawing to be branching. We have a
connected planar map with $f$ faces, each of which is simply connected and has
size at least $3$. (The {\em size} of a face is the number of edges along its
boundary, where an edge is counted twice if both of its sides belong to the
face.) As in the case of simple graphs, we have that $2e$ is equal to the sum
of the sizes of the faces, which is at least $3f$. Hence, by Euler's
polyhedral formula,
$$2=n-e+f\le n-e+{2\over 3}e=n-{1\over 3}e,$$
and the result follows. $\Box$

\medskip

\noindent {\bf Corollary 2.2.}
{\em Let $G$ be a branching topological multigraph with $n\ge 3$ vertices and
$e$ edges. Then for the number of crossings in $G$ we have $c(G)\ge e-3n+6$.}

\smallskip

\noindent {\bf Proof.} By our assumptions, each crossing belongs to precisely
two edges. At each crossing, delete one of these two edges. The remaining
topological graph $G'$ has at least $e-c(G)$ edges. Since $G'$ is a branching
topological multigraph with no two crossing edges, we can apply Lemma 2.1 to
obtain $e-c(G)\le 3n-6.$ \; $\Box$

\medskip

\noindent {\bf Proof of Theorem 3.}
Let $G$ be a branching topological multigraph with $n$ vertices. It is
sufficient to show that for the degree of every vertex $v\in V(G)$ we have
$d(v)\le 2n-4$. This implies that $e(G)\le n(2n-4)/2=n(n-2)$.

Let $v_1, v_2,\ldots, v_{n-1}$ denote the vertices of $G$ different from
$v$. Delete all edges of $G$ that are not incident to $v$. No two remaining
edges cross each other. If $v$ is not adjacent to some $v_i\in V(G)$, then add
a single edge $vv_i$ without creating a crossing. The resulting topological
multigraph, $G'$, is also branching. Starting with any edge connecting $v$ to
$v_1$, list all edges incident to $v$ in clockwise order, and for each edge
write down its endpoint different from $v$. In this way, we obtain a sequence
$\sigma$ of length at least $d(v)$, consisting of the symbols $v_1,
v_2,\ldots, v_{n-1}$, with possible repetition. Let $\sigma'$ denote the
sequence of length at least $d(v)+1$ obtained from $\sigma$ by adding an extra
symbol $v_1$ at the end.

\smallskip

\noindent {\em Property A:} No two consecutive symbols of $\sigma'$ are the same.

This is obvious for all but the last pair of symbols, otherwise the
corresponding pair of edges of $G'$ would form a simple closed Jordan curve
with no vertex in its interior or in its exterior, contradicting the fact that
$G'$ is branching. The last two symbols of $\sigma'$ cannot be the same
either, because this would mean that $\sigma$ starts and ends with $v_1$, and
in the same way we arrive at a contradiction.

\smallskip

\noindent {\em Property B:} $\sigma'$ does not contain a subsequence of the
type $v_i \ldots v_j \ldots v_i \ldots v_j$ for $i\neq j$.

Indeed, otherwise the closed curve formed by the pair of edges connecting $v$
to $v_i$ would cross the closed curve formed by the pair of edges connecting
$v$ to $v_j$, contradicting the fact that $G'$ is crossing-free.

\smallskip

A sequence with Properties A and B is called a {\em Davenport-Schinzel
  sequence of order} $2$. It is known and easy to prove that any such sequence
using $n-1$ distinct symbols has length at most $2n-3$; see~\cite{ShA95}, page
6. Therefore, we have $d(v)+1\le 2n-3$, as required.
\smallskip

To see that the bound in Theorem 3 is tight, place a regular $n$-gon on the
equator $E$ (a great circle of a sphere), and connect any two consecutive
vertices by a single circular arc along $E$. Connect every pair of
nonconsecutive vertices by two half-circles orthogonal to $E$: one in the
Northern hemisphere and one in the Southern hemisphere. The total number of
edges of the resulting drawing is $2{n\choose 2}-n=n(n-2)$. See
Fig. \ref{otszog}. $\Box$

\section{Proof of Theorem 1---using Theorems 2 and 3}

Let $G'$ be a branching topological multigraph of $n'$ vertices and $e'>4n'$ edges.
If $e'\le 10^8n'$, then it follows from Corollary 2.2 that $G'$ meets the requirements of Theorem 1.

To prove Theorem 1, suppose for contradiction that $e'>10^8n'$ and that the number of crossings in $G'$
satisfies
$$c(G')< c(e')^{3}/(n')^{2},$$
for a small constant $c>0$ to be specified later.
\smallskip

Let $d$ denote the average degree of the vertices of $G'$, that is, $d=2e'/n'$. For every
vertex $v\in V(G)$ whose degree, $d(v)$, is larger than $d$, split $v$ into several vertices of
degree at most $d$, as follows. Let $vw_1, vw_2,\ldots, vw_{d(v)}$ be the edges incident
to $v$, listed in clockwise order. Replace $v$ by $\lceil d(v)/d\rceil$ new vertices,
$v_1, v_2,\ldots, v_{\lceil d(v)/d\rceil}$, placed in clockwise order on a very small circle
around $v$. By locally modifying the edges in a small neighborhood of $v$, connect $w_j$ to $v_i$
if and only if $d(i-1)<j\le di$. Obviously, this can be done in such a way that we do not create
any new crossing or two parallel edges that bound a region that contains no vertex. At the end of
the procedure, we obtain a branching topological multigraph $G$ with $e=e'$ edges, and $n<2n'$
vertices, each of degree at most $d=2e'/n'<4e/n$.

Thus, for the number of crossings in $G$, we have
\begin{equation}\label{start}
c(G)=c(G')<4ce^3/n^2
\end{equation}

We break $G$ into smaller components, according to the
following procedure.

\bigskip

\noindent {\sc Decomposition Algorithm}

\medskip

\noindent {\sc Step 0.} {\bf Let} $G^0=G, G^0_1=G, M_0=1, m_0=1.$

Suppose that we have already executed {\sc Step} $i$, and
that the resulting branching topological graph, $G^{i}$, consists of $M_i$ components,
$G_1^{i}, G_2^{i}, \ldots , G_{M_{i}}^{i}$, each having at most
$(4/5)^in$ vertices. Assume without loss of generality
that the first $m_i$ components of $G^i$ have at least
$(4/5)^{i+1}n$ vertices and the remaining $M_i-m_i$ have fewer. Letting
$n(G_j^i)$ denote the number of vertices of the component $G_j^i$, we have
\begin{equation}\label{nulla}
(4/5)^{i+1}n(G)\le n(G_j^{i})\le (4/5)^{i}n(G),\ \ \ \
1\le j\le {m_{i}}.
\end{equation}
Hence,
\begin{equation}\label{darab}
m_i\le (5/4)^{i+1}.
\end{equation}

\smallskip

\noindent {\sc Step $i+1$.} {\bf If}
\begin{equation}\label{egyes}
(4/5)^{i} < {1\over 2}\cdot{e\over n^{2}},
\end{equation}
{\bf then} {\sc stop}. (\ref{egyes}) is called the {\it stopping
rule}.

{\bf Else}, for $j=1, 2, \ldots , {m_{i}}$, delete ${\rm b}_{\rm br}(G_j^{i})$
edges from $G_j^{i}$, as guaranteed by Theorem 2, such that $G_j^{i}$ falls into
two components, each of which is a branching topological graph with at most
$(4/5)n(G_j^{i})$ vertices. Let $G^{i+1}$ denote the resulting topological graph
on the original set of $n$ vertices. Clearly, each component of $G^{i+1}$ has at
most $(4/5)^{i+1}n$ vertices.

\bigskip

Suppose that the {\sc Decomposition Algorithm} terminates
in {\sc Step} $k+1$.
If $k>0$, then
\begin{equation}\label{kettes}
(4/5)^k < {1\over 2}\cdot
{e\over n^{2}} \le (4/5)^{k-1}.
\end{equation}

First, we give an upper bound on the total number of edges
deleted from $G$. Using the fact that, for any nonnegative numbers
$a_1, a_2, \ldots , a_m$,
\begin{equation}\label{harmas}
\sum_{j=1}^m\sqrt{a_j}\le \sqrt{m\sum_{j=1}^m a_j},
\end{equation}
we obtain that, for any $0\le i<k$,
$$\sum_{j=1}^{m_i}\sqrt{c(G_j^{i})}\le
\sqrt{m_i\sum_{j=1}^{m_i}c(G_j^i)}\le
\sqrt{(5/4)^{i+1}}\sqrt{c(G)}<
\sqrt{(5/4)^{i+1}}\sqrt{4ce^{3}/n^{2}}.$$
Here, the last inequality follows from (\ref{start}).

Denoting by $d(v,G^i_j)$ the degree of vertex $v$ in $G^i_j$,
in view of (\ref{harmas}) and (\ref{darab}), we have
$$\sum_{j=1}^{m_i}\sqrt{\sum_{v\in V(G_j^{i})}d^2(v,G^i_j)+n(G^i_j)}\le
\sqrt{m_i\left(\sum_{v\in V(G^i)}d^2(v,G^i)+n\right)}\le$$
$$\sqrt{(5/4)^{i+1}}\sqrt{\max_{v\in V(G^i)}d(v,G^i)\cdot
\sum_{v\in V(G^i)}d(v,G^i)+n}\le
\sqrt{(5/4)^{i+1}}\sqrt{{4e\over n}2e+n}<\sqrt{(5/4)^{i+1}}{3e\over\sqrt{n}}.$$
Thus, by Theorem 2, the total number of edges deleted during the decomposition procedure is

$$\sum_{i=0}^{k-1}\sum_{j=1}^{m_{i}}{\rm b}_{\rm br}(G_j^{i})\le
22\sum_{i=0}^{k-1}\sum_{j=1}^{m_{i}}\sqrt{c(G_j^{i})+\sum_{v\in V(G_j^{i})}d^2(v,G^i_j)+n(G_j^i)}\le$$
$$22\sum_{i=0}^{k-1}\sum_{j=1}^{m_{i}}\sqrt{c(G_j^{i})}
+22\sum_{i=0}^{k-1}\sum_{j=1}^{m_{i}}\sqrt{\sum_{v\in V(G_j^{i})}d^2(v,G^i_j)+n(G_j^i)}\le$$
$$22\left(\sum_{i=0}^{k-1}\sqrt{(5/4)^{i+1}}\right)\left(\sqrt{4ce^3\over n^2}+{3e\over\sqrt{n}}\right)<
350{n\over\sqrt{e}}\left(\sqrt{4ce^3\over n^2}+{3e\over\sqrt{n}}\right)<$$
$$350(2\sqrt{c}e+3\sqrt{en})< 350(2\sqrt{c}e+3\sqrt{e(2e/10^8)})  < {e\over 2},$$
provided that $c\le 10^{-7}$. In the last line, we used our assumption that $e>10^8n'>(10^8/2)n$. The estimate
for the term $\sum_{i=0}^{k-1}\sqrt{(5/4)^{i+1}}$ follows from (\ref{kettes}).
\medskip

So far we have proved that the number of edges of the graph $G^k$ obtained
in the final {step} of the {\sc Decomposition Algorithm} satisfies
\begin{equation}\label{vegso}
e(G^k)> {e\over 2}.
\end{equation}
(Note that this inequality trivially holds if the algorithm terminates
in the very first step, i.e., when $k=0$.)
\smallskip

Next we give a lower bound on $e(G^k)$.
The number of vertices of each connected component of $G^k$ satisfies
$$n(G_j^{k})\le (4/5)^kn < {1\over 2}\cdot
{{e}\over {n^{2}}}n
={e \over 2n},\ \ \ \  1\le j\le M_k.$$
By Theorem 3,
$$e(G_j^{k})\le n^{2}(G_j^{k})<n(G_j^{k})\cdot {e\over 2n}.$$
Therefore, for the total number of edges of $G^k$ we have
$$e(G^k)=\sum_{j=1}^{M_k}e(G_j^{k})<
{e\over 2n}\sum_{j=1}^{M_k}n(G_j^{k})={e\over 2},$$
contradicting (\ref{vegso}). This completes the proof of Theorem 1.
$\Box$

\section{Branching bisection width vs. number of crossings\\
---Proof of Theorem 2}


Suppose that there is a weight function $w$ on a set $V$. Then
for any subset $S$ of $V$, let $w(S)$ denote the
total weight of the elements of  $S$.
We will apply the following separator theorem.

\smallskip

\noindent {\bf Separator Theorem} (Alon-Seymour-Thomas~\cite{AST94}).
 {\em Suppose that a graph $G$ is drawn in the plane with no crossings.
Let $V=\{ v_1, \ldots , v_n\}$ be the vertex set of $G$.
Let $w$ be a nonnegative weight function on $V$.
Then there is a simple closed curve $\Phi$ with the following properties.

(i) $\Phi$ meets $G$ only in vertices.

(ii) $|\Phi \cap V|\le 3\sqrt{n}$

(iii) $\Phi$ divides the plane into two regions, $D_1$ and $D_2$, let
$V_i=D_i\cap V$. Then for $i=1,2$,
$$w(V_i)+{1\over 2}w(\Phi\cap V)\le {2\over 3}w(V).$$ }

\smallskip

Consider a branching drawing of $G$ with exactly $c(G)={\rm cr}_{\rm br}(G)$
crossings.
Let $V_0$ be the set of {\em isolated} vertices of $G$, and let
$v_1, v_2, \ldots , v_m$ be the other vertices of $G$ with
degrees $d_1, d_2, \ldots , d_m$, respectively.
Introduce a new vertex at each crossing. Denote the set of these
vertices by $V_X$.

For $i= 1, 2\ldots , m$,
replace vertex $v_i$
by a set $V_i$ of vertices forming a very small
$d_i\times d_i$ piece of a square
grid, in which each vertex is connected to its
horizontal and vertical neighbors.
Let each edge incident to $v_i$ be hooked up to distinct vertices
along one side of the boundary of $V_i$ without creating any crossing.
These $d_i$ vertices will be called the {\em special boundary vertices}
of $V_i$.

Note that we modified the drawing of the edges only in small neighborhoods of
the grids $V_i$, that is, in nonoverlapping small neighborhoods of
the vertices of $G$, far from any crossing.

\vskip 0.3cm

\begin{figure}[h]
 \begin{center}
  \includegraphics[scale=1.1]{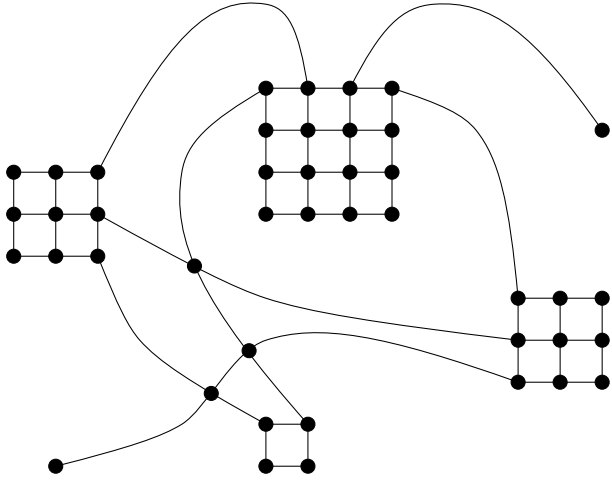}
 \end{center}
\caption{Topological graph $H$.}\label{simonovi1}
\end{figure}


Thus, we obtain a (simple) topological graph $H$, of
$|V_X|+\sum_{i=0}^m|V_i|\le {c}(G)+\sum_{i=1}^md^2_i+n$
vertices and with no crossing; see Fig. \ref{simonovi1}.
For every $1\le i\le m$, assign weight
$1/d_i$ to each special boundary vertex of $V_i$.
Assign weight $1$ to every vertex of $V_0$
and weight $0$ to all other vertices of $H$.
Then $w(V_i)=1$ for every $1\le i\le n$ and $w(v)=1$ for every $v\in V_0$. Consequently,
$w(V(H))=n$.

Apply the Separator Theorem to $H$. Let $\Phi$ denote the closed curve satisfying the conditions of the theorem.
Let $A(\Phi)$ and $B(\Phi)$ denote the region {\em interior} and the {\em exterior} of $\Phi$, respectively.
For $1\le i\le m$, let  $A_i=V_i\cap A(\Phi)$,  $B_i=V_i\cap B(\Phi)$,  $C_i=V_i\cap \Phi$.
Finally, let $C_X=V_X\cap\Phi$.

\smallskip

\noindent{\bf Definition.} For any $1\le i\le m$, we say that

  $V_i$ is of {\em type A} if $w(A_i)\ge {5\over 6}$,

  $V_i$ is of {\em type B} if $w(B_i)\ge {5\over 6}$,

  $V_i$ is of {\em type C,} otherwise.

\noindent For every $v\in V_0$,

  $v$  is of {\em type A} if $v\in A(\Phi)$,

  $v$  is of {\em type B} if $v\in B(\Phi)$,

  $v$  is of {\em type C,} if $v\in\Phi$.
\medskip

Define a partition
$V(G)=V_A\cup V_B$ of
the vertex set of $G$, as follows.
For any $1\le i\le m$, let $v_i\in V_A$ (resp. $v_i\in V_B$)
if $V_i$ is of type $A$ (resp. type $B$).
Similarly,  for every $v\in V_0$, let
$v\in V_A$ (resp. $v\in V_B$)
if $v$ is of type $A$ (resp. type $B$).
The remaining vertices
will be assigned either to $V_A$ or to $V_B$ so as to minimize
$\bigl\lvert|V_A|-|V_B|\bigr\rvert$.

\medskip

\noindent {\bf Claim 4.1}
{\em ${n\over 5}\le |V_A|,|V_B|\le {4n\over 5}$}

\smallskip

\noindent{\bf Proof.} To prove the claim, define another partition
$V(H)=\overline{A}\cup\overline{B}\cup\overline{C}$ such that
$\overline{A}\cap V_i=A\cap V_i$ and $\overline{B}\cap V_i=B\cap V_i$
for $V_0$ and for every $V_i$ of type $C$.
If $V_i$ is of type $A$ (resp. type $B$),
then let $V_i=\overline{A}_i\subset \overline{A}$
 (resp.  $V_i=\overline{B}_i\subset \overline{B}$),
finally, let
$\overline{C}=V(H)-\overline{A}-\overline{B}$.

For any $V_i$ of type $A$, we have $w(\overline{A}_i)-w(A_i)\le {w(A_i)\over 5}$.
Similarly, for any $V_i$ of type $B$, we have $w(\overline{B}_i)-w(B_i)\le {w(B_i)\over 5}$.
Therefore,
$$|w(\overline{A})-w(A)|\le{1\over 5}\cdot\max\{w(A), w(B)\}\le {2n\over 15}.$$
Hence, ${n\over 5}\le w(\overline{A})\le {4n\over 5}$ and,
analogously, ${n\over 5}\le w(\overline{B})\le {4n\over 5}$. In particular,
$|w(\overline{A})-w(\overline{B})|\le {3n\over 5}.$
Using the minimality of
$\bigl\lvert|V_A|-|V_B|\bigr\rvert$, we obtain that
$\bigl\lvert|V_A|-|V_B|\bigr\rvert\le {3n\over 5}$, which implies Claim 4.1. $\Box$

\vskip 0.3cm

\begin{figure}[h]
 \begin{center}
  \includegraphics[scale=0.6]{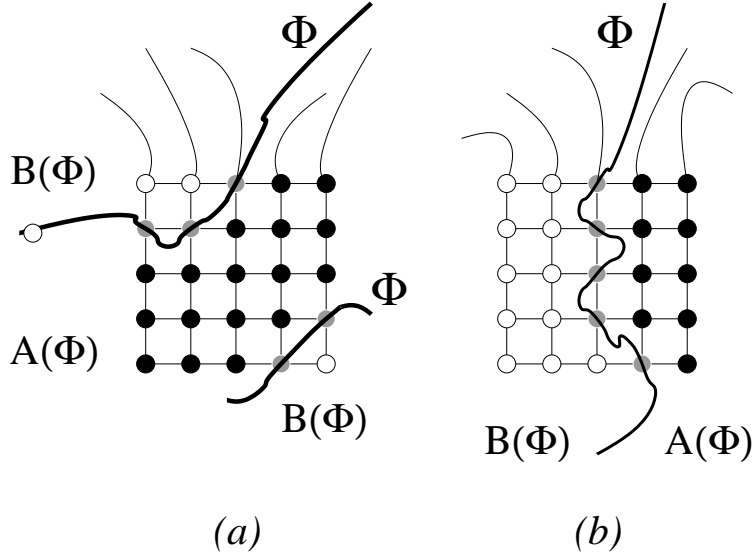}
 \end{center}
\caption{Parts (a) and (b) show a grid of type $A$ and a grid of type $C$, respectively.}\label{typeAC}
\end{figure}

\medskip

\noindent {\bf Claim 4.2.}  {\em For any $1\le i\le n$,

(i)\;  if $V_i$ is of type A (resp. of type B),
then   $|C_i|\ge w(B_i)d_i$
(resp. $|C_i|\ge w(A_i)d_i$);

(ii)\;  if $V_i$ is of type C, then
$|C_i|\ge {d_i\over 6}$.}

\smallskip

\noindent{\bf Proof.} In $V_i$, every connected component belonging to
$A_i$ is separated from every connected component
belonging to $B_i$ by vertices in $C_i$.
There are $w(A_i)d_i$ (resp. $w(B_i)d_i$) special boundary vertices
in $V_i,$ which belong to $A_i$ (resp. $B_i$).
It can be shown by an easy case analysis that
the number of separating points
$|C_i|\ge \min\{w(A_i), w(B_i)\}d_i$, and Claim 4.2 follows; see
Fig. \ref{typeAC}. $\Box$

\medskip

\noindent {\bf Claim 4.3.}
{\em Let $V=V(G)$. There is a closed curve
$\Psi$, not passing through any vertex of $H$, whose interior and exterior are denoted by $A(\Psi)$ and $B(\Psi)$, resp., such that

 (i)\;\;\; $V\cap A(\Psi)=V_A$,

 (ii)\;\;  $V\cap B(\Psi)=V_B$,

 (iii)\; the total number of edges of $G$ intersected by $\Psi$ is at most
$$18\sqrt{{c}(G)+\sum_{i=1}^nd_i^2+n}.$$
}

\smallskip

\noindent{\bf Proof.}
For any $1\le i\le m$, we say that

$V_i$ is of {\em type 1}\;\;  if $|C_i|\ge d_i/6$, 

$V_i$ is of {\em type 2}\;\; if $|C_i|<d_i/6$.

\noindent For every $v\in V_0$,

 $v$ is of {\em type 1}\;\; if $v\in \Phi$,

 $v$ is of {\em type 2}\;\; if $v\in A(\Phi)\cup B(\Phi)$.

\smallskip

\noindent It follows from Claim 4.2 that if a set $V_i$ or an isolated vertex $v\in V_0$
is of type C, then it is also of type 1.

Next, we modify the curve $\Phi$
in small neighborhoods of the grids $V_i$ and of the isolated vertices $v\in V_0$ to make sure
that the resulting curve $\Psi$  satisfies the conditions in the claim.


Assume for simplicity
that $v_i\in V_A$; the case $v_i\in V_B$ can be treated analogously.
If $v_i$ is a vertex of degree at most $1$ and $\Phi$
passes through $v_i$, slightly perturb $\Phi$ in a small neighborhood of $v_i$
(or slightly shift $v_i$) so that after this change $v_i$ lies in the interior
of $\Phi$.
Suppose next that
the degree of $v_i$ is at least $2$.
Let $S_i$ and
$S'_i\subset S_i$ be two closed squares containing $V_i$ in their interiors,
and assume that $S_i$ (and, hence, $S'_i$) is only slightly larger than the
convex hull of the vertices of $V_i$.
We
distinguish two cases.


\smallskip

\vskip 0.3cm

\begin{figure}[h]
 \begin{center}
\includegraphics[scale=0.6]{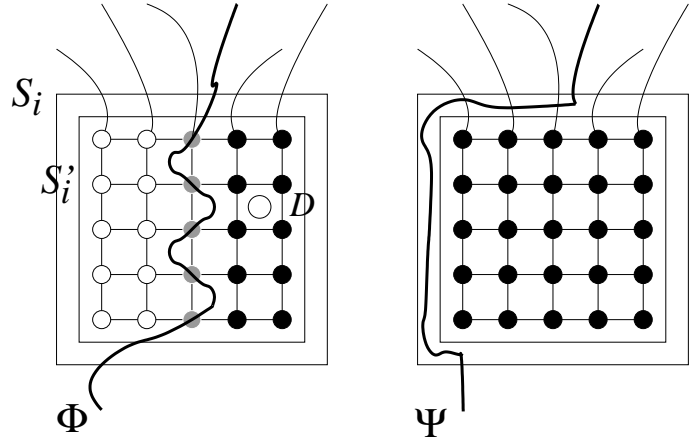}
 \end{center}
\caption{Claim 4.3, Case 1.}
\end{figure}



{\sc Case 1.}  {\em $V_i$ is of type 1.}
Let $D$ be a small disk in $S'_i$ that belongs to the interior of
$\Phi$ and let $p$ be its center.
Let $\tau : S_i\rightarrow S_i$ be a
homeomorphism of $S_i$ to itself which
keeps the boundary of $S_i$ fixed and let $\tau(D)=S'_i$.
%
%
%
Observe that every piece of $\Phi$ within the convex hull
of the vertices of $V_i$ is mapped into an arc in the very narrow ring
$S_i\setminus S'_i$. In particular, if we keep the vertices and the edges of the
grid $H[V_i]$ (as well as all other parts of the drawing) fixed, after this
local modification $\Phi$ will avoid all vertices of $V_i$ and it may
intersect only those (at most $d_i$) edges incident to $V_i$ which correspond
to original edges of $G$ and end at some special boundary vertex of
$V_i$. Moreover, after this modification, every vertex of $V_i$ will lie in $A(\Phi)$,
in the {\em interior} of $\Phi$.
\smallskip

\vskip 0.3cm

\begin{figure}[h]
 \begin{center}
\includegraphics[scale=0.6]{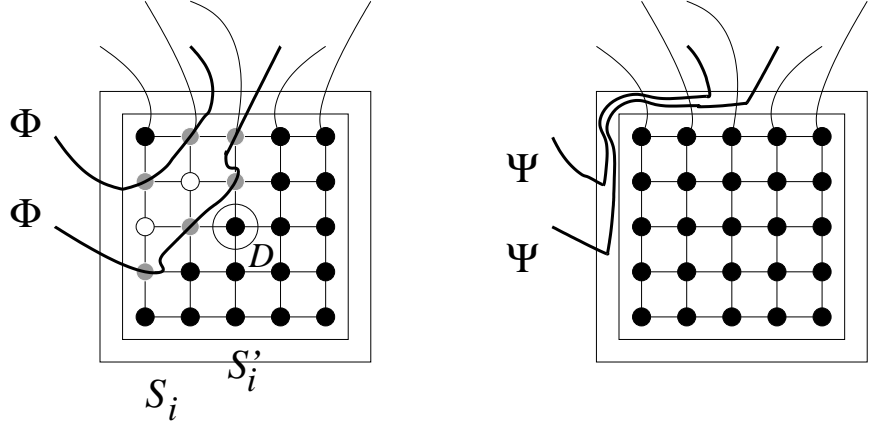}
 \end{center}
\caption{Claim 4.3, Case 2.}
\end{figure}

{\sc Case 2.} {\em $V_i$ is of type 2.}
In this case, by Claim 4.2, $V_i$ is of type $A$.

Orient $\Phi$ arbitrarily. Let $(p_1, p_1'), (p_2, p_2'), \ldots $ denote
the point pairs at which $\Phi$ enters and leaves the convex hull of $V_i$, so that the arc between $p_jp'_j$
lies inside the convex hull of $V_i$, for every $j$. Note that both $p_j$ and $p'_j$ are vertices of $V_i$.
In view of the fact that $|C_i|\le d_i/6$, we know that the (graph) distance between $p_j$
and $p_j'$ (in $H[V_i]$) is at most $d_i/6$.
More precisely, for every $j$, the points $p_j$ and $p_j'$ divide the boundary of the convex hull of $V_i$
into two arcs. We call the shorter of these arcs
the {\em boundary interval
defined by $p_j$ and $p_j'$}, and denote it by $[p_j, p_j']$.
By assumption, the {\em length} of $[p_j, p_j']$. the number of edges of $H[V_i]$ comprising $[p_j, p_j']$,  is at most $d_i/6$.

It is not hard to see that the curve $\Phi$ cannot came close to the center
$p$ of $V_i$ and that $p$ belongs to the interior of $\Phi$.
Let $D$ be a small disk centered at $p$. Then $D$
also
belongs to the interior of
$\Phi$.
Let $\tau : S_i\rightarrow S_i$ be a
homeomorphism of $S_i$ to itself such that
(i) $\tau$ keeps the boundary of $S_i$ fixed, (ii) $\tau(D)=S'_i$,
(iii) $\tau(p)=p$,
and
(iv) for any $q\in S_i$, that points  $p$, $q$, and $\tau(q)$ are collinear.
Observe that every piece  $(p_j, p_j')$,
of $\Phi$ within the convex hull
of the vertices of $V_i$
is mapped into an arc in the very narrow ring
$S_i\setminus S'_i$, along the corresponding
boundary interval, $[p_j, p_j']$,
defined by $p_j$ and $p_j'$.
In particular, if we keep the vertices and edges of the
grid $H[V_i]$ (as well as all other parts of the drawing) fixed, after this
local modification $\Phi$ will avoid all vertices of $V_i$ and it may
intersect only those (at most $d_i/6$) edges incident to $V_i$ which correspond
to original edges of $G$ and end at some special boundary vertex of
$V_i$ in a boundary interval. Moreover, now every vertex of $V_i$ will lie {\em inside} $\Phi$.

\smallskip

Repeat the above local modification for each $V_i$ and for each $v\in V_0$.
The resulting curve, $\Psi$, satisfies conditions (i) and (ii).
It remains to show that it also satisfies (iii).







\smallskip



To see this, denote by $E_X$ the set of all edges of $H$ adjacent
to at least one element of $C_X$.
For any $1\le i\le m$, define $E_i\subset E(H)$
as follows.
If $V_i$ is of type 1, then let all edges of $H$ leaving
$V_i$ belong to $E_i$.
If $V_i$ is of type 2, then by Claim 4.2,
it can be of type A or B, but not C.
Let $E_i$ consist of
all edges leaving
$V_i$ and crossed by $\Psi$.

For any $1\le i\le m$, let  $E'_i$ denote
the set of edges of $G$ corresponding to
the elements of $E_i$ ($0\le i\le m$) and let
$E'_X$ denote
the set of edges corresponding to
the elements of $E_X$.

Clearly, we have $|E'_i|\le |E_i|,$ because distinct edges of
$G$ give rise to distinct edges of $H$.
Since $V_A$ and $V_B$ are on different sides of $\Psi$, it crosses all edges  between
$V_A$ and $V_B$.

Obviously, $|E'_X|\le |E_X|\le 4|C_X|$.
By Claim 4.2, if $V_i$ is of type 1,
then $|E'_i|\le |E_i|=d_i\le 6|C_i|$.
If $V_i$ is of type 2,
then $|E'_i|\le |E_i|=d_i\le |C_i|$.
Therefore,
$$|E(V_A, V_B)|\le |\cup_{i=0}^nE'_i|\le \sum_{i=0}^n|E_i|\le
6|C|\le 18\sqrt{{c}(G)+\sum_{i=1}^nd_i^2+n}.$$

This finishes the proof of Claim 4.3. $\Box$

\medskip

Now we are in a position to complete the proof of Theorem 2.
Remove those edges of $G$ that are cut by $\Psi$. Let
$G_A$ (resp. $G_B$)
be the subgraph of the resulting graph $G'$, induced by $V_A$ (resp. $V_B$), with the inherited drawing.
Suppose that, e.g., $G_B$ is not a branching topological graph.
Then it has an {\em empty lens}, that is, a region bounded by two parallel edges
that does not contain any vertex of $V_B$. There are two types of empty lenses: bounded and unbounded.
We show that there are at most $\sqrt{{c}(G)}$ bounded empty lenses, and at most
$\sqrt{{c}(G)}$ unbounded empty lenses in $G_B$.

Suppose that $e$ and $e'$ are two parallel edges between $v$ and $v'$ which enclose a bounded empty lens $L$.
Then $v$ and $v'$ are in the exterior of $\Psi$, and $\Psi$ does not cross the edges $e$ and $e'$.
As $G$ was a branching topological multigraph, both $L$ and its complement contain at least one vertex of $G$
in their interiors. Since $L$ is empty in $G_B$, it follows that all vertices of $G$ inside $L$ must belong to $V_A$,
and, hence, must lie in the interior of $\Psi$. Thus, $\Psi$ must lie entirely inside the lens $L$.

Suppose now that $f$ and $f'$ are two
other parallel edges between two vertices $u$ and $u'$, and they determine
another bounded empty lens $M$. Arguing as above, we obtain
that $\Psi$ must also lie entirely inside $M$.
Then $v$ and $v'$ are outside of $M$, and $u$ and $u'$ are outside of $L$.
Therefore, 
these four edges determine four crossings. Any such
crossing can belong to only one pair of bounded empty lenses $\{L,M\}$,
we conclude that for the number of
bounded empty lenses $k$ in $G_B$ we have
$4{f\choose 2}\le c(G)$, therefore,
$k\le\sqrt{{c}(G)}$.
Analogously, there are at most $\sqrt{{c}(G)}$ unbounded empty lenses in $G_B$.

We can argue in exactly the same way for $G_A$. Thus, altogether there are at most
$4\sqrt{{c}(G)}$ empty lenses in $G_A$ and $G_B$.
If we delete a boundary edge of each of them, then no empty lens is left.

Thus, by deleting the edges of $G$ crossed by $\Psi$ and then one boundary edge of each empty lens,
we obtain a decomposition of $G$ into two branching topological multigraphs, and the
number of deleted edges is at most

$$18\sqrt{{c}(G)+\sum_{i=1}^nd_i^2+n}+4\sqrt{{c}(G)}\le
22\sqrt{{c}(G)+\sum_{i=1}^nd_i^2+n}.$$

This concludes the proof of Theorem 2. $\Box$

\medskip

{\bf Acknowledgement.} We are very grateful to Stefan Felsner,
Michael Kaufmann, Vincenzo Roselli,
Torsten Ueckerdt, and Pavel Valtr for their many valuable comments,
suggestions, and for many interesting discussions during the Dagstuhl Seminar
"Beyond-Planar Graphs: Algorithmics and Combinatorics", November 6-11, 2016,
\\
http://www.dagstuhl.de/en/program/calendar/semhp/?semnr=16452.

\end{document}